\documentclass[<options>]{elsarticle}

\usepackage{CJK}

\usepackage{amsfonts}
\usepackage{amsthm}
\usepackage{latexsym}
\usepackage{enumerate}
\usepackage{graphicx}
\usepackage{amssymb,amsmath}

\usepackage{float}
\usepackage{texdraw}
\usepackage[dvips]{psfrag}
\usepackage{amsmath}
\usepackage{indentfirst}

\usepackage{titlesec}
\titleformat{\section}{\Large}{\thesection.}{0.4em}{\textbf}{}
\titleformat{\subsection}{\large}{\thesubsection}{0.5em}{\textbf}{}
\usepackage{lineno,hyperref}
\journal{Journal of \LaTeX\ Templates}
\modulolinenumbers[5]

\newtheorem{Lemma}{Lemma}[section]

\newtheorem{Definition}{Definition}[section]

\newcommand\R{{\bf R}}

\begin{document}

\begin{frontmatter}

\title{\textbf{\LARGE{On the number of hyperelliptic limit cycles of Li\'{e}nard systems}}}
%\tnotetext[mytitlenote]{Fully documented templates are available in the elsarticle package on \href{http://www.ctan.org/tex-archive/macros/latex/contrib/elsarticle}{CTAN}.}

%% Group authors per affiliation:
\author{\textbf{Xinjie Qian }}
\cortext[]{Corresponding author.}
\ead{1701110014@pku.edu.cn}
\author{ \textbf{Jiazhong Yang}}
\ead{~~jyang@math.pku.edu.cn}
\address{School of Mathematical Science, Peking University, 100871 Beijing,  P. R. China}

%\author{\textbf{XinJie Qian , JiaZhong Yang}}

%\address {School of Mathematical Science, Peking University, Beijing, P. R. China}

%\maketitle

%\rule[-10pt]{16.3cm}{0.05em}
%\renewcommand{\abstractname}

\begin{abstract}

In this paper, we study the maximum number, denoted by $H(m,n)$, of  hyperelliptic limit cycles
 of the Li\'{e}nard systems
$$\dot x=y, \qquad \dot y=-f_m(x)y-g_n(x),$$
where, respectively, $f_m(x)$ and $g_n(x)$ are real polynomials of degree $m$  and $n$, $g_n(0)=0$.
The main results of the paper are as follows: We obtain the upper bound and lower bound  of $H(m,n)$ in all the  cases with  $n\neq 2m+1$. When $n=2m+1$, we  derive the lower bound of $H(m,n)$. Furthermore, these upper bound can be  reached in some cases.
\end{abstract}
%\vspace{6mm}

%\noindent{\it Keywords\/}: Hyperelliptic Limit Cycles, Li\'{e}nard Systems

%\rule[-10pt]{16.3cm}{0.05em}


\end{frontmatter}

%\linenumbers

\section {Introduction}
Consider the  following   Li\'{e}nard  differential system
\begin{equation}\label{1}
\dot x=y, \quad\quad\dot y=-f_m(x)y-g_n(x),
\end{equation}
where $f_m(x)$ and $ g_n(x)$ are polynomials of degrees $m$ and  $n$, respectively, with the following explicit expressions
   $$f_m(x)=\sum_{i=0}^ma_ix^i, \qquad  \ g_n(x)=\sum_{i=1}^nb_ix^i, \quad a_mb_n\not=0.$$
We shall call this system  a {\it Li\'{e}nard   system of type $(m,n)$},
or simply   a Li\'{e}nard  system if no confusion arises.

This paper is primarily devoted to a study
of the  maximum number $H(m,n)$ of hyperelliptic limit cycles of the Li\'enard system in terms of $m$ and $n$.

Here we adopt the conventional definition of a limit cycle. Namely,
by {\it a limit cycle}   of a polynomial system we mean that it is
  an isolated closed orbit of the system.
It is called an {\it algebraic limit cycle} if it is a limit cycle and  is contained in
an invariant algebraic curve \{$(x,y)\mid F(x,y)=0$\}.
In particular, if  $F(x,y)$ takes the form
$F(x,y)=(y+P(x))^2-Q(x),$
where $P$ and $Q$ are polynomials,
then we call the invariant curve
hyperelliptic.
Correspondingly,  a limit cycle is called a  hyperelliptic limit cycle if it is contained in a hyperelliptic curve.

The  investigation   of  limit cycles of the Li\'enard system has been one
of the most interesting topics for decades (see \cite{JL},\cite{FD}).
In the most general setting, however, it is a very hard subject
 and the  problem of existence  is quite
elusive. Therefore certain assumptions are reasonably imposed, and  special categories are
technically restricted.  Among them,    the algebraic and hyperelliptic versions of the problem have caught
particular attention of the study.  A brief survey of the situation is as follows.

Odani \cite{Odani} in 1995  proved that if $n\leq m$ and
$f_mg_n(f_m/g_n)'\not\equiv0$,
then  any  Li\'{e}nard system of  $(m,n)$-type  has   no  invariant algebraic curves.
Therefore in this case, it is  impossible to have
any hyperelliptic limit cycles.

Chavarriga et al. \cite{JIJH}, Zoladek \cite{Z}, and Makoto Hayashi \cite {MH}  proved that
any  Li\'{e}nard  systems of  the  types
 $(0,n)$,   $(1,n)$, $(2, 4)$  and  $(m,m +1)$
 have  no algebraic limit cycles, hence there are no hyperelliptic limit cycles.

In 2008, Llibre and Zhang  \cite{LZ} proved that   no Li\'enard
system of   $(3,5)$-type  has  hyperelliptic limit cycles.
On the other hand, in the same paper \cite{LZ}, they   found
 that in the following cases there are  Li\'{e}nard systems   of $(m,n)$-type  which
can possess  hyperelliptic limit cycles:
\begin{itemize}

\item [(i)] $(m, n)$-type, for    $m\geq 2$ and $n\geq 2m+1$;
\item[(ii)] $(m,2m)$-type for  $m\geq 3$;
\item[(iii)] $(m, 2m-1)$-type for $m \geq 4$;
\item[(iv)]   $(m, 2m-2)$-type for $m\geq 4$.
\end{itemize}

An individual type $(5,7)$
of the Li\'enard system  is discussed  in \cite{YZ},
where    Yu and Zhang  clarified that
there exist   Li\'{e}nard systems  of   $(5,7)$-type
which  have hyperelliptic limit cycles.

A recent paper \cite{Liu} is conclusive, where the authors considered  the remaining types of the systems and proved that, in all these cases, there  always exist  Li\'enard  systems
  of $(m,n)$-type which   have  hyperelliptic limit cycles.
Thus the problem of the existence of   hyperelliptic limit cycles
for all the possible types of the Li\'enard systems is
completely answered.

Collecting all the known results mentioned above and arranging them into Fig.1,
 we   can provide a visual
way to exhibit the distribution of the  hyperelliptic limit cycles. Namely,  in the $(m,n)$-plane,
there is a clearly-cut boundary
dividing all the types of the Li\'enard systems into two regions:
Systems falling in region $1$ can
never have any  hyperelliptic limit cycle which means $H(m,n)=0$, and
in the other region, for each pair of
$(m,n)$,  there always exists   such  a Li\'enard system
which admits  at least one   hyperelliptic limit cycle, thus $H(m,n)\geq1$.
Systems falling on the boundary are also unambiguously  specified.

\begin{figure}
  \centering
  \includegraphics[width=2.5in]{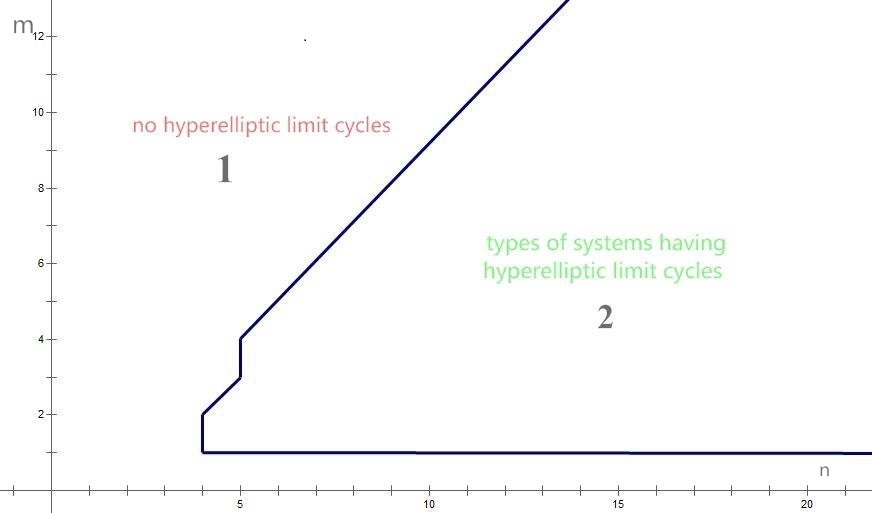}\\
\label{fig}{Fig.1.The maximum number of hyperelliptic limit cycles.}
\end{figure}

The present paper grows from a very casual observation.
If one looks at Figure and  takes region $1$ as land and region $2$ as sea, and if
we walk from the land to the sea,
we are  in fact traveling  from a region  where systems have  no hyperelliptic limit cycle
to a region where  such limit cycles start to appear.
A very  natural question like this can  pop up:
when we walk from the land to the sea,
does the water become deeper and deeper? In other words,
does the maximum number of hyperelliptic limit cycles
increase as we walk into  the sea  further and further?
Such curiosity  leads us to  explore  this problem and
to  see if there is any  algebraic mechanism  behind this.
The investigation  turns    out to be quite interesting:
%Indeed, this plain feeling is  quite soundable!
While only   those Li\'enard  systems falling in the sea
 can  have hyperelliptic
limit cycles, we prove that
 those systems in ``deeper" water indeed can have larger $H(m,n)$.
A detailed classification is summarized in the following theorem.
Notice that we also consider the configuration of these limit cycles,
another one of very important aspects of the subject.

\vspace{4mm}
{\bf \large{Main Theorem}}: {\it
Consider  Li\'enard systems of the type $(m,n)$ where $m\geq 2$, the maximum number of hyperelliptic limit cycles admits the following estimations:

 $$ H(m,n)\geq\left\{
\begin{array}{lcl}
 n-m-1       &      & {m+2 \leq n \leq [\frac{4m+2}{3}]}\\
{[\frac{n-1}{4}]} &      & {[\frac{4m+2}{3}]+1\leq n\leq 2m, m\geq 4}\\
{[\frac{m}{2}]}   &      & {n\geq 2m+1}
\end{array} \right.
$$

and
$$ H(m,n)\leq\left\{
\begin{array}{lcl}
{[\frac{n+1}{4}]}     &      & {m+2 \leq n \leq 2m-2, m\geq 4}\\
{[\frac{n-1}{4}]}  &     & {n=2m-1, or~n=2m , m\geq 4}\\

{[\frac{m}{2}]}  &      &   {n>2m+1}
\end{array} \right. $$

In all the cases  with  $ n\neq {2m+1}$ and $H(m,n)>1$,
 the hyperelliptic limit cycles of  the  system  can  only have non-nested configuration.
}

\vspace{4mm}
{\bf \large{ Remark: }}{\it  It immediately follows from the main theorem that
\begin{itemize}
\item [\emph{(i)}] When $1+[\frac{4m+2}{3}]\leq n\leq 2m-2$, if $n-1\equiv0$ \emph{(mod $4$)} or $n-1\equiv1$ \emph{(mod $4$)}, then $H(m,n)=[\frac{n-1}{4}]$;
  \item[\emph{(ii)}] If $n=2m-1$ or $n=2m$, $m\geq4$, then $H(m,n)=[\frac{n-1}{4}]$;
  \item[\emph{(iii)}] If $ n>{2m+1}$ , then $H(m,n)=[\frac{m}{2}]$.

\end{itemize}
}

 The paper is organized as follows: In section 2, we shall introduce
some preliminaries including definitions, notation
and basic methods. In section 3 ,4 and 5, we present a detailed proof of the results.

\section{Preliminaries}

In this section, we  shall collect some related properties of Li\'enard
systems and introduce a complete discrimination system for polynomials. For the proof of these results, we refer the reader to
(\emph{\emph{\cite{ Liu, LZ, YZ, Z, MH, YL, YLJZ}}})for details.

\subsection{Hyperelliptic limit cycles of Li\'enard systems}
Recall that the Li\'enard  system  takes the form
$$\dot x=y, \quad\quad\dot y=-f_m(x)y-g_n(x).$$
Assume that the system   has a hyperelliptic invariant curve
\begin{equation}\label{2}
F(x, y)= (y + P(x))^2-Q(x)=0.
\end{equation}
The following properties hold, whose proof is standard and hence omitted.

\begin{Lemma}
 There exists $K(x, y)\in \R[x, y]$ such that
$$y\frac{\partial F}{\partial x}-(f_m(x)y+g_n(x))\frac{\partial F}{\partial
y}=K(x,y)F.$$

\end{Lemma}

\begin{Lemma} If  relation \emph{(\ref{2})} holds, then
the degree of  polynomial $P(x)$ is $m+1$, and   the polynomials  $f_m$ and $g_n$
 can be expressed in terms of $P$ and $Q$ as follows.
 \begin{equation}\label{3}
f_m=P'+\frac{PQ'}{2Q}, \quad
g_n=\frac{Q'(P^2-Q)}{2Q}.\end{equation}
\end{Lemma}

  Since  any
 singular point   of system (\ref{1}) must be   located on the $x$-axis, thus
when  a hyperelliptic curve $F(x, y)=0$ contains a limit cycle of system (\ref{1}),
the  limit cycle   should intersect the $x$-axis at two different points, denoted
by $(s_1, 0)$ and $(s_2, 0)$. The following properties hold.

\begin{Lemma}

\emph{(i)} $s_1$ and $s_2$ are real  simple roots of $Q (x)$.\emph{ (ii)}   Any  root of $Q(x)$ must be   a root of
$P(x)$.

 \end{Lemma}

 %Moreover, with the above assumption, the following conclusion is obvious.

\begin{Lemma}\label{lemma2}
If $s_1$ and $s_2$ are simple roots of $Q (x)$ and
$Q(x)>0$ in $(s_1, s_2)$, then the hyperelliptic curve \emph{(\ref{2})}
contains a closed curve in the strip $s_1\leq x\leq s_2$.

 \end{Lemma}

Now one step further: assume that
(i) the hyperelliptic curve $F(x,y)=0$ contains a
closed curve $C$ in the strip $s_1\leq x\leq s_2$, where $s_1$ and $s_2$
are simple roots of $Q (x)$,
(ii) this closed curve $C$ surrounds only
one singularity $(\alpha, 0)$ of system  (\ref{1}),
(iii) the
singularity $(\alpha, 0)$ is a focus or node.
Then  this closed
curve $C$ is a limit cycle.

We have the following criteria  to recognize the type of the singular point.
 \begin{Lemma}
 If $g_n(\alpha)=0, g_n'(\alpha)>0$ and $f_m(\alpha)\not=0$, then
$(\alpha, 0)$ is  a focus or a node of system  \emph{(\ref{1})} .
 \end{Lemma}

Combining all the known result, we give the following lemma which is very useful in determining if an algebraic curve
 is a hyperelliptic limit cycle of the Lienard system.
  \begin{Lemma} \label{lemma4}
An algebraic   curve \emph{(\ref{2})} in the strip $x\in [s_1,s_2]$
is a hyperelliptic limit cycle if the following sufficient conditions
are met.

\emph{(i)}  $f_m$ and $g_n$  satisfy (\ref{3}),

\emph{(ii)} All the roots of $Q(x)=0$ are real and  $s_1,s_2$ are simple root
and $Q(x)>0$ for $x\in (s_1,s_2)$.

\emph{(iii)}  $P^2(x)-Q(x)<0$ for
$x\in (s_1,s_2)$.

\emph{(iv) } If $\alpha \in (s_1,s_2)$ such that $Q'(\alpha)=0$, then $f_m(\alpha)\neq 0$.

  \end{Lemma}

Proof of Lemma:  Condition (i) means that  $F=0$ is the invariant curve
of the system, and all the roots of $Q(x)$ are the roots of $P(x)$.
From (iii) we know that  the curve $F(x,y)=0$ in the strip
bounded by $x=s_1$  and $x=s_2$  intersects the $x$ axis only at these two
endpoints.  Condition (ii) means that $Q'(s_i)\neq 0$. It follows that the curve in the strip
has no singular points and is closed.  From (ii)
we also know that   $Q'(x)$ has only one real root $\alpha$ for $[s_1,s_2]$.
Again,  from (iii) we see that  $g_n(x)$ has a unique real root $\alpha$.
Therefore the system has only one singular point
inside the closed orbit formed by $F=0$ when restricted to the strip.
We can even see that this singular  point is either a focus type or
a node. In fact, $g_n'(\alpha)=Q''(\alpha)\cdot \frac{P^2(\alpha)-Q(\alpha)}{2Q(\alpha)}
+Q'(\alpha)(\frac{P^2-Q}{2Q})'(\alpha)$.
Notice that the second term vanishes, and since $\alpha$ is the maximal
value point of $Q$, therefore $Q''(\alpha)<0$. It follows that
$g'(\alpha)>0$. Condition (iv) says that
$f_m(\alpha)\neq 0$, therefore $(\alpha,0)$ is
a focus or a node. Therefore the closed orbit
is hyperelliptic limit cycle of the system.

\subsection{Algorithm for root classification}

 Given a polynomial
      $$f(x)=a_0x^n+a_1x^{n-1}+\cdots+a_n,$$
 we write the derivative of $f(x)$ as
 $$f'(x)=na_0x^{n-1}+(n-1)a_1x^{n-2}+\cdots+a_{n-1}.$$

For the n-degree polynomial $f(x)$, $\alpha_1, \alpha_2, \cdots, \alpha_n$ denote all the roots of it. Let $s_p=\sum\limits_{j=1}^n\alpha_j^p$, $p=0,1,2,\cdots,n$,  $S_k=|s_{i+j}|,i,j=0,1,\cdots,k-1,$ that is,
\begin{equation}
S_k=\left|
 \begin{array}{cccc}
 s_0 & s_1 & \cdots & s_{k-1} \\
 s_1 &  s_2 & \cdots & s_{k} \\
 \cdots & \cdots & \cdots & \cdots \\
 s_{k-1} & s_{k} & \cdots & s_{2k-2}\\
 \end{array}
\right|.
\end{equation}
\vskip0.2cm
\begin{Definition}
\emph{(discrimination matrix)}  The Sylvester matrix of $f(x)$ and $f'(x)$,  denoted by \emph{Discr(f)}
$$
\left(
 \begin{array}{cccccccc}
a_0 & a_1 & a_2 & \cdots & a_n & 0 & \cdots & 0 \\
0 & na_0 & (n-1)a_1& \cdots & a_{n-1} & 0 & \cdots & 0  \\
0 & a_0 & a_1 & \cdots & a_{n-1} & a_n & \cdots & 0 \\
 0 & 0 & na_0 & \cdots & 2a_{n-2} & a_{n-1} & \cdots & 0 \\
\cdots & \cdots & \cdots & \cdots & \cdots & \cdots & \cdots & \cdots \\
\cdots & \cdots & \cdots & \cdots & \cdots & \cdots & \cdots & \cdots \\
0 & 0 & 0 & \cdots & a_0 & a_1 & \cdots & a_n\\
0& 0 & 0 & \cdots & 0 & na_0 & \cdots & a_{n-1}\\
 \end{array}
 \right)
$$
is called the discrimination matrix of $f(x)$.
\end{Definition}
\vskip0.2cm
\begin{Definition}
\emph{(discriminant sequence)}  Denoted by $D_k$, the determinant of the submatrix of \emph{Discr(f)}, formed by the first 2k rows and the first $2k$ columns, for $k=1,\cdots,n$. We call the n-tuple \emph{$(D_1,D_2,\cdots,D_n)$}
 the discriminant sequence of polynomial $f(x)$.
 \end{Definition}
\vskip0.2cm
\begin{Definition}
 \emph{(sign list)}
we call the list
$$\emph{$[ sign(D_1), sign(D_2), \cdots, sign(D_n)]$}$$
the sign list of the discrimination sequence \emph{$(D_1,D_2,\cdots,D_n)$}.
\end{Definition}
 \vskip0.2cm
 \begin{Definition}
\emph{(revised sign list)} Given a sign list $[s_1,s_2,\cdots,s_n]$, we construct
a new list $[\varepsilon_1,\varepsilon_2,\cdots,\varepsilon_n]$ as follows:
\begin{itemize}
\item If $[s_1,s_2,\cdots,s_n]$ is a section of given list, where
$s_i\neq0, s_{i+1}=\cdots=s_{i+j-1}=0, s_{i+j}\neq0$, then we replace the subsection $[s_{i+1},s_{i+2},\cdots,s_{i+j-1}]$
     by $[-s_i,-s_i,s_i,s_i,-s_i,-s_i,s_i,s_i,-s_i,\cdots]$.

 i.e. let
$\varepsilon_{i+r}=(-1)^{[{\frac{r+1}{2}}]}s_i$, for $r=1,2,\cdots{j-1}.$

\item Otherwise, let $\varepsilon_k=s_k$, there are no changes for other terms.
\end{itemize}
\end{Definition}
 \vskip0.2cm
 From \emph{\cite{YLJZ}}, we already know the following lemma.
\vskip0.35cm
\begin{Lemma} \label{lemma5}
 For $k=1,2,\cdots,n,$ we have $D_k=S_k.$
   \end{Lemma}
 \vskip0.35cm  \begin{Lemma} \label{lemma6}

 Given a polynomial \emph{$f(x)=a_0x^n+a_1x^{n-1}+\cdots+a_n$} with real coefficients, if the number of the sign changes of the revised sign list of
$$\{D_1(f),D_2(f),\cdots,D_n(f)\}$$
is $v$, then the number of the pairs of distinct conjugate imaginary roots of \emph{$f(x)$ } equals $v$.
Furthermore, if the number of non-vanishing members of the revised sign list is $l$, then the number of the distinct real roots of \emph{$f(x)$} equals $l-2v$.
  \end{Lemma}

\section{The Proof of the Results about Lower Bounds}
According to all the possible pairs $(m,n)$ where $m\geq 2$, we divide the proof into the following cases.
\begin{itemize}
\item[] (i)  $ m+2\leq n\leq [\frac{4m+2}{3}]$;
\item[] (ii) $[\frac{4m+2}{3}]+1\leq n\leq 2m$ and $(m,n)$ is not in $\{(3,5),(2,4)\}$;
\item[] (iii)$n\geq 2m+1$.
\end{itemize}

\subsection{Case (i)}
When $ m+2\leq n\leq [\frac{4m+2}{3}]$, it suffices to construct
a Li\'enard system of type $(m,n)$ which can have $n-m-1$ hyperelliptic limit
cycles on invariant curve ({\ref 2}).

Suppose $n$ is odd. Now let $t=\frac{4m-3n+3}{2}$. By Corollary 3.1 in
\cite{Liu}, there exist a positive constant $c$ and a polynomial
\begin{displaymath}
Q_1(x)=(x-x_0)(x-1)\prod^{t}_{i=1}(x-x_i)^2\prod^{n-m-2}_{i=1}(x-y_i)^2,
\end{displaymath}
such that
\begin{displaymath}
P_1(x)=Q_1(x)+c=\prod^{t}_{i=1}(x-z_i)^2\prod^{n-m-1}_{i=1}(x-a_i)(x-b_i),
\end{displaymath}
where $x_0<z_1<x_1<z_2<...<z_t<x_t$ and $x_t<a_1<b_1<y_1<a_2<b_2<y_2<...<a_{n-m-1}<b_{n-m-1}<1$. We set

$$G(x)=(x-x_0)^2(x-1)^2\prod^{t}_{i=1}(x-x_i)^2\prod^{n-m-2}_{i=1}(x-y_i)^2\prod^{n-m-1}_{i=1}(x-a_i)(x-b_i),$$
\begin{displaymath}
\begin{aligned}
P(x)&=\sqrt{G(x)P_1(x)}\\
 &=(x-x_0)(x-1)\prod^{t}_{i=1}(x-x_i)(x-z_i)\prod^{n-m-2}_{i=1}(x-y_i)\prod^{n-m-1}_{i=1}(x-a_i)(x-b_i),
\end{aligned}
\end{displaymath}
\begin{displaymath}
\begin{aligned}
Q(x)&=G(x)Q_1(x)\\
&=(x-x_0)^3(x-1)^3\prod^{t}_{i=1}(x-x_i)^4\prod^{n-m-2}_{i=1}(x-y_i)^4\prod^{n-m-1}_{i=1}(x-a_i)(x-b_i).
\end{aligned}
\end{displaymath}%We take $f_{m}(x)$ and $g_{n}(x)$ in system ({\ref 1}) in the form of equation ({\ref 3}),
then
$$
f_{m}(x)=P'(x)+\frac{P(x)Q'(x)}{2Q(x)},\quad g_{n}(x)=\frac{Q'(x)(P^2(x)-Q(x))}{2Q(x)}
$$
are polynomials of degree $m$ and $n$ respectively.

We claim, for each $i$, $i=1,2,...,n-m-1$, when $x\in[a_i,b_i]$, the closed curve given by ({\ref 2})
is a hyperelliptic limit cycle of the system.

1. In fact, it is easy to see that the
condition(i), (ii), (iii) of Lemma \ref{lemma4} is satisfied.

2. Let us verify condition(iv) by contradiction. Assume $Q'(x)$ and $f_{m}(x)$ have a common root $\alpha$ in $(a_i,b_i)$, then $P'(\alpha)=0$. With $G(\alpha)\neq0$, then $G'(\alpha)=P_1'(\alpha)=Q_1'(\alpha)=0$, and $(\frac{G(x)}{P_1(x)Q_1(x)})'\Big| _{x=\alpha}=0$.
\begin{displaymath}
\left(\frac{G(x)}{P_1(x)Q_1(x)}\right)'=\frac{(x-x_0)(x-1)}{\prod_{i=1}^{t}(x-z_i)^2}\left(\frac{1}{x-x_0}+\frac{1}{x-1}-2\sum_{i=1}^{t}\frac{1}{x-z_i}\right),
\end{displaymath}
we have $(G/(P_1Q_1))'(\alpha)>0$
, this leads to a contradiction.

By Lemma \ref{lemma4}, we can prove the system has
$n-m-1$ hyperelliptic limit cycles.

Suppose $n$ is even, let $t=(4m-3n+2)/2$. By Corollary 3.2 in
\cite{Liu}, there exist a positive constant $c$ and a polynomial
\begin{displaymath}
Q_1(x)=(x-1)\prod^{t}_{i=1}(x-x_i)^2\prod^{n-m-1}_{i=1}(x-y_i)^2,
\end{displaymath}
such that
\begin{displaymath}
P_1(x)=Q_1(x)+c=(x-x_0)\prod^{t}_{i=1}(x-z_i)^2\prod^{n-m-1}_{i=1}(x-a_i)(x-b_i),
\end{displaymath}
where $x_0<x_1<z_1<x_2<...<x_t<z_t$ and $z_t<y_1<a_1<b_1<y_2<a_2<b_2<...<y_{n-m-1}<a_{n-m-1}<b_{n-m-1}<1$. We set
$$G(x)=(x-x_0)(x-1)^2\prod^{t}_{i=1}(x-x_i)^2\prod^{n-m-1}_{i=1}(x-y_i)^2\prod^{n-m-1}_{i=1}(x-a_i)(x-b_i),$$
\begin{displaymath}
\begin{aligned}
P(x)&=\sqrt{G(x)P_1(x)}\\
&=(x-x_0)(x-1)\prod^{t}_{i=1}(x-x_i)(x-z_i)\prod^{n-m-1}_{i=1}(x-y_i)\prod^{n-m-1}_{i=1}(x-a_i)(x-b_i),
 \end{aligned}
\end{displaymath}
\begin{displaymath}
\begin{aligned}
Q(x)&=G(x)Q_1(x)\\
&=(x-x_0)(x-1)^3\prod^{t}_{i=1}(x-x_i)^4\prod^{n-m-1}_{i=1}(x-y_i)^4\prod^{n-m-1}_{i=1}(x-a_i)(x-b_i).
\end{aligned}
\end{displaymath}
then
$$
f_{m}(x)=P'(x)+\frac{P(x)Q'(x)}{2Q(x)},\quad g_{n}(x)=\frac{Q'(x)(P^2(x)-Q(x))}{2Q(x)}
$$
are polynomials of degree $m$ and $n$ respectively.

We claim, for each $i$, $i=1,2,...,n-m-1$, when $x\in[a_i,b_i]$, the closed curve given by ({\ref 2})
is a hyperelliptic limit cycle of the system.

1. In fact, it is easy to see that the
condition(i), (ii), (iii) of Lemma \ref{lemma4} is satisfied.

2. Let us verify condition(iv) by contradiction.
Assume $Q'(x)$ and $f_{m}(x)$ have a common root $\alpha$ in $(a_i,b_i)$. Analogous the argument
above, we can get $\alpha$ is a root of $(G/(P_1Q_1))'$, while
\begin{displaymath}
\left(\frac{G(x)}{P_1(x)Q_1(x)}\right)'=\frac{(x-1)}{\prod_{i=1}^{t}(x-z_i)^2}\left(\frac{1}{x-1}-2\sum_{i=1}^{t}\frac{1}{x-z_i}\right),
\end{displaymath}
for each $i$, we can observe $x_0<z_i<\alpha<1$, thus $(G/P_1Q_1)'(\alpha)>0$ ,
this leads to a contradiction.

By Lemma \ref{lemma4}, we can prove the system has $n-m-1$ hyperelliptic limit cycles.

\subsection{Case(ii)}

Now we come to case \emph{$(ii)$}, when $[\frac{4m+2}{3}]+1\leq n\leq 2m-1$,
we shall  construct
a Li\'enard system ({\ref1}) that can have $[\frac{n-1}{4}]$ hyperelliptic limit
cycles on invariant curve ({\ref 2}), from which we can infer that $H(m,n)\geq [\frac{n-1}{4}]$.

In the proof of case $(i)$, we perturbed the polynomial with a constant to transform
repeated roots into single roots, but this perturbation doesn't work in case $(ii)$. To prove case $(ii)$,
firstly we divide the case $(ii)$ into the following cases:
\begin{itemize}
\item[] (ii-i)  $n-1\equiv0$ (mod $4$);
\item[] (ii-ii) $n-1\equiv1$ (mod $4$);
\item[] (ii-iii) $n-1\equiv2$ or $n-1\equiv3$ (mod $4$);
\end{itemize}
{\bf Case (ii-i)}: $n-1\equiv0$ (mod $4$)
\begin{Lemma}\label{lemma7}
For $h,l\in N$, define the polynomial
\begin{displaymath}
Q_1(x)=(x-s)x^{2h+1}\prod_{i=1}^l(x-i)^2,
\end{displaymath}
where $s>l+1$, then  there exists a polynomial $c(x)$ of degree $2h$ which is positive in $[0,s]$ and such that
\begin{displaymath}
Q_1(x)+c(x)=(x-y_{l+1})\prod_{i=1}^l(x-y_i)(x-z_i)\prod_{i=1}^{2h+1}(x-x_i)
\end{displaymath}
where $0<x_1<x_2<\cdots<x_{2h+1}<y_1,$ $ y_1<1<z_1<y_2<\cdots<z_l<y_{l+1}<s$.
\end{Lemma}
{\bf Proof:}
We prove this lemma by mathematical induction. For
$h=0$, let $c(x)$ be a positive constant $\epsilon$. It easily follows
that the proposition for $h=0$ holds, if $\epsilon$ is sufficiently
small. Assume the proposition holds for $h=k$, it must been shown
that the proposition holds for $h=k+1$.

Decompose $Q_1(x)$
into two fractions $x^2$ and $Q_1^*(x)$, then $0$ is a repeated
root of degree of $2k+1$ of $Q_1^*(x)$. Using the induction
hypothesis, there exists a polynomial $c^*(x)$ of degree $2k$
which is positive in $[0,s]$, (we can choose $c^*(x)$ which satisfied the maximum absolute value of its coefficients is sufficient small) and such that
\begin{eqnarray}
Q_1(x)+x^2c^*(x)&=&x^2(Q_1^*(x)+c^*(x)) \nonumber\\
                &=&x^2(x-y_{l+1}')
                    \prod_{i=1}^l(x-y_i')(x-z_i')\prod_{i=3}^{2k+3}(x-x_i'), \nonumber
\end{eqnarray}
where  $0<x_3'<\cdots<x_{2k+3}'<y_1'$, $y_1'<1<z_1'<y_2'<\cdots<z_l'<y_{l+1}'<s$.

Choose a sufficiently small $d$ which satisfied  $d>0$ and $xc^*(x)-d$ has only one root $\alpha<<1$ in $[0,s]$. For the maximum absolute value of coefficients of $c^*(x)$ is sufficient small, the local maximum of $Q_1(x)+x^2c^*(x)$ in
$(0,x_3')$ is the least maximum among all the maxima of
$Q_1(x)+x^2c^*(x)$ in $[0,s]$. Perturbing $Q_1(x)+x^2c^*(x)$ with $-dx$,
we get a polynomial
\begin{displaymath}
Q_1(x)+x^2c^*(x)-dx= x(x-y_{l+1}'')
                    \prod_{i=1}^l(x-y_i'')(x-z_i'')\prod_{i=2}^{2k+3}(x-x_i''),
\end{displaymath}
where $0<x_2''<x_3''<\cdots<x_{2k+3}''<y_1''$, $y_1''<z_1''<y_2''<\cdots<z_l''<y_{l+1}''$.

Since $\alpha<<1$ is the only root of  $xc^*(x)-d$ in $[0,s]$, we have $Q_1(s)+s^2c^*(s)-ds>0$ and $Q_1(i)+i^2c^*(i)-di>0, 1\leq i\leq l$ , then $y_{l+1}''<s$ and $y_i''<i<z_i'', 1\leq i\leq l$. When $0<x<\alpha$, we have $x^2c^*(x)-dx<0$, while $x_2''$ is the root of $Q_1(x)+x^2c^*(x)-dx$, for $Q_1( x_2'')<0$, then $\alpha< x_2''$.

Assume $\gamma$ is minimum point of $x^2c^*(x)-dx$  in $[0,s]$, then $0<\gamma<\alpha<x_2''$. Choose $b>0$ satisfied $\gamma^2c^*(\gamma)-d\gamma+b>0$, $Q_1(\gamma)+\gamma^2c^*(\gamma)-d\gamma+b<0$. (The existence of $b$ relies on $Q_1(\gamma)<0$.) Now we start to proof all roots of $Q_1(x)+x^2c^*(x)-dx+b$ are real. Assume $\beta$ is minimum point of $Q_1(x)+x^2c^*(x)-dx+b$ in $[0,x_2'']$, we obtain $Q_1(\beta)+\beta^2c^*(\beta)-d\beta+b\leq Q_1(\gamma)+\gamma^2c^*(\gamma)-d\gamma+b<0$. For $d$ is sufficiently small, the local minimum $Q_1(\beta)+\beta^2c^*(\beta)-d\beta+b$ in
$(0,x_2'')$ is the largest minimum among all the minima of
$Q_1(x)+x^2c^*(x)-dx$ in $[0,s]$, %and $Q_1(x^*)+{x^*}^2c^*(x^*)-dx^*+b>0$ where $x^*$ is the root of $Q_1(x)+x^2c^*(x)-dx$,
we know  that
all roots of $Q_1(x)+x^2c^*(x)-dx+b$ are real.

Perturbing $Q_1(x)+x^2c^*(x)-dx$ with $b$, we get a polynomial
\begin{displaymath}
Q_1(x)+x^2c^*(x)-dx+b=(x-y_{l+1})\prod_{i=1}^l(x-y_i)(x-z_i)\prod_{i=1}^{2k+3}(x-x_i).
\end{displaymath}
where $0<x_1<x_2<\cdots<x_{2k+3}<y_1$,  $y_1<y_1''<z_1''<z_1<y_2<\cdots<z_l<y_{l+1}<y_{l+1}''$. For $y_{l+1}''<s$ and $y_i''<i<z_i'', 1\leq i\leq l$, we have $0<x_1<x_2<\cdots<x_{2k+3}<y_1$,$y_1<1<z_1<y_2<\cdots<z_l<y_{l+1}<s$. On the other hand, we know  the degree of $c(x)=x^2c^*(x)-dx+b$ is $2k+2$, and $c(x)\geq \gamma^2c^*(\gamma)-d\gamma+b>0$ in $[0,s]$. This completes the proof of lemma.

\vskip 0.4cm
Denote $\frac{n-1}{4}=t$, We set
\begin{displaymath}
Q_1(x)=(x-2m+2t)x^{6t-2m+1}\prod_{i=1}^{m-2t-1}(x-i)^2,
\end{displaymath}

by Lemma \ref{lemma7}, we can perturb $Q_1(x)$ with a polynomial $c(x)$ of degree $6t-2m$ which is positive in $[0,s]$, then
\begin{displaymath}
P_1(x)=Q_1(x)+c(x)= \prod_{i=1}^{t}(x-a_i)(x-b_i),
\end{displaymath}
where
$0<a_1<b_1<\cdots<a_{3t-m+1}$, $a_{3t-m+1}<b_{3t-m+1}<1<a_{3t-m+2}<\cdots <m-2t+1<a_t<b_t<2m-2t$. We define

$$G(x)=\prod_{i=1}^t(x-a_i)(x-b_i)\prod_{i=0}^{m-2t-1}(x-i)^2(x-2m+2t)^2.$$
\begin{displaymath}
P(x)=\sqrt{G(x)P_1(x)},\qquad Q(x)=G(x)Q_1(x),
\end{displaymath}
then
$$P(x)=(x-2m+2t)\prod_{i=1}^{t}(x-a_i)(x-b_i)\prod_{i=0}^{m-2t-1}(x-i),$$
\begin{displaymath}
Q(x)=(x-2m+2t)^3x^{6t-2m+3}\prod_{i=1}^{t}(x-a_i)(x-b_i)\prod_{i=1}^{m-2t-1}(x-i)^4,
\end{displaymath}
and

$$f_{m}(x)=P'(x)+\frac{P(x)Q'(x)}{2Q(x)},\quad
g_{n}(x)=\frac{Q'(x)(P^2(x)-Q(x))}{2Q(x)}$$
are polynomials of degree $m$ and $n$ respectively.

 We claim, for each $i$, $i=1,2,...,t$, when $x\in[a_i,b_i]$,
 the closed curve given by ({\ref 2}) is a hyperelliptic limit cycle of the system.

 1. In fact, it is easy to see that the condition(i), (ii), (iii) of Lemma \ref{lemma4}
 is satisfied.

 2. Let us verify condition(iv) by contradiction. Assume $Q'(x)$
 and $f_{m}(x)$ have a common root $\alpha$ in $(a_i,b_i)$.

Suppose $6t-2m=0$, then $n=\frac{4m+3}{3}$, which is possible when $\frac{4m+3}{3}=[\frac{4m+2}{3}]+1$.

 For  $G(\alpha)\neq0$, we get $\alpha$
 is the common root of $P_1'$, $Q_1'$, $G'$ and $(G/(P_1Q_1))'$, then

 \begin{displaymath}
\left(\frac{G}{P_1Q_1}\right)'=[x^2(x-2m+2t)/x^{6t-2m+1}]'=2x-2m+2t,
\end{displaymath}
thus $\alpha=m-t$, but it is impossible for $Q_1'(m-t)<0$ and this leads
to a contradiction.

By Lemma \ref{lemma4}, we can prove the system has $t$
hyperelliptic limit cycles.

On the other hands, $6t-2m>0$, for $6t-2m$ is even, then $6t-2m\geq 2$. With $f_{m}(\alpha)=Q'(\alpha)=0$, then $P'(\alpha)=0$. For $G(\alpha)\neq0$, and $\frac {P_1}{Q_1}=\frac {P^2}{Q}$, we get $(\frac {Q}{P})'\Big|_{x=\alpha}=(\frac {P_1}{Q_1})'\Big|_{x=\alpha}=0$. Since
\begin{displaymath}
\frac{Q(x)}{P(x)}=x^{6t-2m+2}(x-2m+2t)^2\prod_{i=1}^{m-2t-1}(x-i)^3,
\end{displaymath}
we know $\alpha$ is irrelevant of $c(x)$.
Differentiating $P_1/Q_1$, we have
\begin{eqnarray}\label{5}
\left(\frac{P_1(x)}{Q_1(x)}\right)'=\frac{c'(x)Q_1(x)-Q_1'(x)c(x)}{Q_1^2(x)}
\end{eqnarray}
For $Q_1(\alpha)\neq0$,
we have \begin{eqnarray}\label{6}\frac{c'(\alpha)}{c(\alpha)}=\frac{Q_1'(\alpha)}{Q_1(\alpha)}.\end{eqnarray}
 With the degree of $c(x)$ is more than 2 and the right side of (\ref{6}) is irrelevant of $c(x)$ , we can change the polynomial coefficients of $c(x)$ to make the left hand side of (\ref{6}) doesn't equal the right hand side, such that the root of (\ref{5})
in $(a_j, b_j)$ is different to the root of equation $(Q/P)'$. Therefore, such $\alpha$ doesn't exist and
this verifies condition(iv).

By Lemma \ref{lemma4}, we prove the system has $t$
hyperelliptic limit cycles. This completes the proof of the case  $n-1\equiv0$ (mod $4$).
\vskip 0.4cm
\noindent{\bf Case (ii-ii)}: $n-1\equiv1$ (mod $4$)
\vskip 0.2cm
\noindent For the proof of Lemma \ref{lemma8} is similar to Lemma \ref{lemma7}, we omit it.
\begin{Lemma}\label{lemma8}
For $h\in N^+, l\in N$, define the polynomial
\begin{displaymath}
Q_1(x)=(x-s_1)(x-s_2)x^{2h}\prod_{i=1}^l(x-i)^2,
\end{displaymath}
where $ s_1<-1, s_2>l+1$, then there exists a polynomial $c(x)$ of degree $2h-1$ which is positive in $[s_1,s_2]$ and such that
\begin{displaymath}
Q_1(x)+c(x)=(x-z_{-1})(x-y_{l+1})\prod_{i=1}^l(x-y_i)(x-z_i)\prod_{i=1}^{2h}(x-x_i)
\end{displaymath}
where $s_1<z_{-1}<x_1<0<x_2<x_3<\cdots<x_{2h}<y_1$, $y_1<1<z_1<y_2<\cdots<z_l<y_{l+1}<s_2$.
\end{Lemma}
Denote $\frac{n-2}{4}=t$.
We set
\begin{displaymath}
Q_1(x)=(x+2)(x-s)x^{6t-2m+2}\prod_{i=1}^{m-2t-2}(x-i)^2,
\end{displaymath}where $s>>m-2t-2$.
Since $6t-2m+2\geq2$, by lemma \ref{lemma8}, we can perturb $Q_1(x)$ with a polynomial $c(x)$ of degree $6t-2m+1$ which is positive in $[-2,s]$, then
\begin{displaymath}
P_1(x)=Q_1(x)+c(x)= \prod_{i=1}^{t}(x-a_i)(x-b_i),
\end{displaymath}
where
$-2<a_1<b_1<0<a_2<\cdots<a_{3t-m+2}<b_{3t-m+2}$, $b_{3t-m+2}<1<a_{3t-m+3}<b_{3t-m+3}<2<a_{3t-m+4}<\cdots <m-2t-2<a_t<b_t<s$. Define
\begin{displaymath}
G(x)=\prod_{i=1}^t(x-a_i)(x-b_i)\prod_{i=0}^{m-2t-2}(x-i)^2(x+2)^2(x-s)^2,
\end{displaymath}
\begin{displaymath}
P(x)=\sqrt{G(x)P_1(x)},\qquad Q(x)=G(x)Q_1(x),
\end{displaymath}
we have
\begin{displaymath}
P(x)=(x+2)(x-s)\prod_{i=1}^{t}(x-a_i)(x-b_i)\prod_{i=0}^{m-2t-2}(x-i),
\end{displaymath}
\begin{displaymath}
Q(x)=(x+2)^3(x-s)^3x^{6t-2m+4}\prod_{i=1}^{t}(x-a_i)(x-b_i)\prod_{i=1}^{m-2t-2}(x-i)^4.
\end{displaymath}
then

$$f_{m}(x)=P'(x)+\frac{P(x)Q'(x)}{2Q(x)},\quad
g_{n}(x)=\frac{Q'(x)(P^2(x)-Q(x))}{2Q(x)}$$
are polynomials of degree $m$ and $n$ respectively.

We claim, for each $i$, $i=1,2,...,t$, when $x\in[a_i,b_i]$,
 the closed curve given by ({\ref 2}) is a hyperelliptic limit cycle of the system.

1. In fact, it is easy to see that the condition(i), (ii), (iii) of Lemma \ref{lemma4}
 is satisfied.

2. Let us verify condition(iv) by contradiction. Assume $Q'(x)$
 and $f_{m}(x)$ have a common root $\alpha$ in $(a_i,b_i)$, then $P'(\alpha)=0$.
 Note that $\frac {P_1}{Q_1}=\frac {P^2}{Q}$, and $G(\alpha)\neq0$, we get $(\frac {Q}{P})'\Big|_{x=\alpha}=(\frac {P_1}{Q_1})'\Big|_{x=\alpha}=0$. Since
\begin{displaymath}
\frac{Q(x)}{P(x)}=x^{6t-2m+3}(x+2)^2(x-s)^2\prod_{i=1}^{m-2t-2}(x-i)^3,
\end{displaymath}
we get $\alpha$ is irrelevant of $c(x)$ immediately.
Differentiating $P_1/Q_1$, then
\begin{eqnarray}\label{7}
\left(\frac{P_1(x)}{Q_1(x)}\right)'=\frac{c'(x)Q_1(x)-Q_1'(x)c(x)}{Q_1^2(x)},
\end{eqnarray}
for $Q_1(\alpha)\neq0$,
we have \begin{eqnarray}\label{8}\frac{c'(\alpha)}{c(\alpha)}=\frac{Q_1'(\alpha)}{Q_1(\alpha)}.\end{eqnarray}
 While $6t-2m+1>0$, the degree of $c(x)$ is more than 1 and the right side of (\ref{8}) is irrelevant of $c(x)$, we can change the polynomial coefficients of $c(x)$ to make the left hand side of (\ref{8}) doesn't equal the right hand side, such that the root of (\ref{7})
in $(a_j, b_j)$ is different to the root of equation $(Q/P)'$.
Therefore, such $\alpha$ doesn't exist and
this verifies condition(iv).

By Lemma \ref{lemma4}, we can prove the system has $t$
hyperelliptic limit cycles. Since $t=\frac{n-2}{4}=[\frac{n-1}{4}]$, we complete the proof.

\vskip 0.4cm
\noindent{\bf Case (ii-iii)}: $n-1\equiv2$ or $n-1\equiv3$(mod $4$)
\begin{Lemma}\label{lemma9}
If a Li\'enard system of $(m,n)$-type has $t$ hyperelliptic limit cycles on
invariant curve $(y + P(x))^2 - Q(x)=0$ and for each limit cycle the conditions
of Lemma \ref{lemma4} are met, then there exists Li\'enard system of $(m+1,n+2)$-type
with at least $t$ hyperelliptic limit cycles.
\end{Lemma}
{\bf Proof.} It suffices that, based on the system in the assumption, we construct
a new Li\'enard system of $(m+1,n+2)$-type in the form of
\begin{displaymath}
\dot { x } = y, \qquad \dot { y } = -f_{m+1}(x)y - g_{n+2}(x)
\end{displaymath}
with the same number of hyperelliptic limit cycles. We take $\tilde P_s(x)$
and $\tilde Q_s(x)$ in the form
\begin{displaymath}
\tilde P_s(x)=P(x)(x-s),\quad \tilde Q_s(x)=Q(x)(x-s)^2.
\end{displaymath}
Changing $P(x)$ and $Q(x)$ in equation (\ref{3}) to $\tilde P_s(x)$ and
$\tilde Q_s(x)$ respectively, we get
\begin{displaymath}
f_{m+1}(x)=\tilde P_s'(x)+\frac{\tilde P_s(x)\tilde Q_s'(x)}{2\tilde Q_s(x)},\quad g_{n+2}(x)=\frac{\tilde Q_s'(x)(\tilde P_s^2(x)-\tilde Q_s(x))}{2\tilde Q_s(x)}.
\end{displaymath}
Note they are polynomials of $m+1$, $n+2$ degree respectively.

Consider a hyperelliptic limit cycle of the original system on the
invariant curve (\ref{2}) that intersect with $x$-axis on $a_1$ and $b_1$. We claim there exists a sufficient large $s_1$,  which satisfied the closed curve with $x\in [a_1,b_1]$ on invariant curve $(y +\tilde P_{s_1}(x))^2-\tilde Q_{s_1}(x)=0$  is a hyperelliptic
limit cycle of the new system.

We observe that the condition(i),
(ii) and (iii) of Lemma \ref{lemma4} are trivially verified when $s$ is larger
than all the roots of $Q(x)$. Then we just have to consider
condition(iv). Differentiating $\tilde Q_s(x)$, we get
\begin{displaymath}
\tilde Q_s'(x)=(x-s)^2(Q'(x)+\frac{2}{x-s}Q(x)).
\end{displaymath}
It follows that, $\tilde \alpha_s \rightarrow \alpha$ as $s\rightarrow\infty$,
where $\alpha$ and $\tilde \alpha_s $ denote the root of $Q'(x)$ and
$\tilde Q_s'(x)$ in $(a_1,b_1)$ respectively.
Differentiating $\tilde P_s(x)$, we get
\begin{displaymath}
\tilde P_s'(x)=(x-s)(P'(x)+\frac{1}{x-s}P(x)).
\end{displaymath}
Hence, $\tilde P_s'(\tilde \alpha_s)/(\tilde \alpha_s-s)\rightarrow P'(\tilde \alpha_s)\rightarrow P'(\alpha)$
as $s\rightarrow\infty$. Furthermore,  $P'(\alpha)\neq 0$ which follows from the assumption that
$f_m'(\alpha)\neq 0$. Thus, we can find a sufficient large $s_1$ satisfied %$\tilde P'(\tilde \alpha_{s_1})/(\tilde
$\tilde P_{s_1}'(\tilde \alpha_{s_1})\neq 0$
to make $f_{m+1}(\tilde \alpha_{s_1})\neq 0$ . By Lemma \ref{lemma4}, we can prove the system $$\dot { x } = y,\quad \dot { y } = -f_{m+1}(x)y - g_{n+2}(x)$$ has at least $t$ hyperelliptic limit cycles ,where \begin{displaymath}
f_{m+1}(x)=\tilde P_{s_1}'(x)+\frac{\tilde P_{s_1}(x)\tilde Q_{s_1}'(x)}{2\tilde Q_{s_1}(x)},\quad g_{n+2}(x)=\frac{\tilde Q_{s_1}'(x)(\tilde P_{s_1}^2(x)-\tilde Q_{s_1}(x))}{2\tilde Q_{s_1}(x)},
\end{displaymath}this completes the proof of the lemma.

\vskip0.2cm
{\bf Proof of the case $n-1\equiv2$ (mod $4$)\emph{}}:
Suppose $(m-1,n-2)$ is still in case $(ii)$, then  $(m-1,n-2)$ is in the case (ii-i) , use the above argument, we have a Li\'enard system of $(m-1,n-2)$-type that has $[\frac{n-3}{4}]$
hyperelliptic limit cycles. Since $n-1\equiv2$ (mod $4$), we have $[\frac{n-3}{4}]=[\frac{n-1}{4}]$. By the argument of Lemma \ref{lemma9} , we can construct a new Li\'enard system of
$(m,n)$-type with  $[\frac{n-1}{4}]$ hyperelliptic limit cycle based on the system of $(m-1,n-2)$-type.

On the other hand, $(m-1,n-2)$ is in case $(i)$, then $[\frac {4m+4}{3}]=[\frac {4m+5}{3}]=n$, but $n-1\equiv2$ (mod $4$), which yields a contradiction. Therefore $(m-1,n-2)$ can only in case $(ii)$, this completes the proof.

\vskip0.2cm
{\bf Proof of the case  $n-1\equiv3$ (mod $4$))}:
Suppose $(m-1,n-2)$ is still in case $(ii)$, then  $(m-1,n-2)$ is in the case (ii-ii), use the above argument, we have a Li\'enard system of $(m-1,n-2)$-type that has $[\frac{n-3}{4}]$
hyperelliptic limit cycles. Since $n-1\equiv3$ (mod $4$), we have $[\frac{n-3}{4}]=[\frac{n-1}{4}]$. By the argument of Lemma \ref{lemma9} , we can construct a new Li\'enard system of
$(m,n)$-type with at least $[\frac{n-1}{4}]$ hyperelliptic limit cycles based on the system of $(m-1,n-2)$-type, thus $H(m,n)\geq[\frac{n-1}{4}]$.

On the other hand, $(m-1,n-2)$ is in case $(i)$, we can construct a Li\'enard system of $(m-1,n-2)$ type that has $n-m-2$
hyperelliptic limit cycles. For $(m-1,n-2)$ is in case $(i)$, and  $n-1\equiv3$ (mod $4$), we have $n-m-2=[\frac{n-1}{4}]$. This completes the proof of the case  $n-1\equiv3$ (mod $4$).

\vskip0.3cm
When $n=2m$,
  we define
 $$P(x)=\prod\limits_{i=1}^m(x-i)(x+s)\quad Q(x)=\prod\limits_{i=1}^m(x-i)(x+s)^{m+2},$$
where $s>>m$ is sufficiently large. If $m$ is odd, for each $i=1,2,\cdots\frac{m-1}{2}$, when $x\in[2i-1,2i]$, the closed curve
given by (\ref{2}) is a hyperelliptic limit cycle of the system, therefore $H(m,2m)\geq\frac{m-1}{2}=[\frac{2m-1}{4}]$.

On the other hand,  $m$ is even, for each  $i=1,2,\cdots\frac{m-2}{2}$, when $x\in[2i,2i+1]$, the closed curve
given by (\ref{2}) is a hyperelliptic limit cycle of the system, therefore $H(m,2m)\geq\frac{m-2}{2}=[\frac{2m-1}{4}]$.

\subsection{Case(iii)}

We set
 $$P(x)=\prod_{i=1}^m(x-i)(x+s),$$
$$Q(x)=-s\prod_{i=1}^m(x-i)(x+s)^{n-m+1},$$ where $s>>m$ is sufficiently large, we take $f_{m}(x)$ and $g_{n}(x)$ in system ({\ref 1}) in the form of equation ({\ref 3}).
It is easy to see $f_m(x)$ and $g_n(x)$ are polynomials of degree $m$ and $n$ respectively.

Suppose $m$ is even. We claim, for each $i=1,2,...\frac{m}{2}$, when $x\in[2i-1,2i]$, the closed curve
given by (\ref{2}) is a hyperelliptic limit cycle of the system.

1. In fact, it is easy to see
that the condition(i), (ii), (iii) of Lemma \ref{lemma4} is satisfied.

2. Let us verify condition(iv)
by contradiction. Assume $Q'(x)$ and $f_{m}(x)$ have a common root $\alpha$ in $(2i-1,2i)$, then $P'(\alpha)=0$. With

\begin{displaymath}
 R'(x)=(\frac{Q(x)}{P(x)})'=-s(n-m)(x+s)^{n-m-1},
\end{displaymath}
we would have $R'(\alpha)=0$, but $R'(x)$ only have one root $-s$, this leads to a contradiction.

By Lemma \ref{lemma4}, we can prove the system
has $\frac{m}{2}$ hyperelliptic limit cycles.

Suppose $m$ is odd. In an analogous way, when $x\in[2i,2i+1]$, we can prove the closed curve
given by (\ref{2}) is a hyperelliptic limit cycle of the system  for each
$i=1,2,\cdots\frac{m-1}{2}$.

Therefore, we obtain $H(m,n)\geq \left[\frac{m}{2}\right]$, when $m\geq2$ and $n\geq2m+1$.
\section{Configuration Of  Hyperelliptic Limit Cycles}

\begin{Lemma}\label{lemma10}
If  an $(m,n)$-Lienard system \emph{(\ref{1})}
has a hyperelliptic curve
$$(y+P(x))^2-Q(x)=0,$$
where $n\neq2m+1$,
then the system  only  has this one hyperelliptic curve.
\end{Lemma}

{\bf Proof.} From equation (\ref{3}), we have
  \begin{equation}\label{9}
2Q(x)f_m(x)=2Q(x)P'(x)+P(x)Q'(x),
\end{equation}
and
\begin{equation}\label{10}
2Q(x)g_n(x)=Q'(x)(P^2(x)-Q(x)).
\end{equation}
Therefore, we know that  the degree of $P(x)$ is $m+1$, while the degree
 of $P^2(x)-Q(x)$ is $n+1$.
Let $f_m(x)$ and $g_n(x)$ take the form
\begin{displaymath}
f_m(x)=\sum_{i=0}^{m}a_i x^i,\qquad g_n(x)=\sum_{i=0}^{n}b_i x^i.
\end{displaymath}

If $n>2m+1$, the degree of $P^2(x)-Q(x)$ equals the degree of $Q(x)$.
Let us denote $P(x)$ and $Q(x)$ by
\begin{displaymath}
P(x)=\sum_{i=0}^{m+1}p_i x^i,\qquad Q(x)=\sum_{i=0}^{n+1}q_i x^i.
\end{displaymath}
then the coefficients of the highest degree terms of the each side of equations (\ref{9}) and (\ref{10}) are:
 \begin{displaymath}
2a_m q_{n+1}=(2m+n+3)p_{m+1}q_{n+1},
\quad 2q_{n+1}b_n=-(n+1)q_{n+1}^2.
\end{displaymath}
Thus
$p_{m+1}=\frac{2a_m}{2m+n+3}$ and $q_{n+1}=\frac{-2b_n}{n+1}$ are uniquely
determined.

Comparing  the coefficients of the second highest degree terms of the polynomials on each side of equations
(\ref{9}) and (\ref{10}) respectively, we have
\begin{displaymath}
2a_m q_n+2a_{m-1}q_{n+1}=(2m+2+n)p_{m+1}q_n+(2m+n+1) p_m q_{n+1},
\end{displaymath}
\begin{displaymath}
2b_n q_n+ 2b_{n-1}q_{n+1}=-(2n+1) q_n q_{n+1}+(n+1)q_{n+1}c,
\end{displaymath}
where  $c=0$ or $c=p_{m+1}^2$. For the coefficients of $p_m$ and $q_n$ in the linear equations mentioned above which  derive from comparing the coefficients of the second highest degree terms of the polynomials on each side of equations
(\ref{9}) and (\ref{10}) are $(2m+n+1)q_{n+1}$ and $nq_{n+1}$ respectively, we have  the values of $p_m$ and $q_n$ are uniquely defined.

More generally, by comparing the coefficients of $x^{n+i}$ and $x^{2n+i-m}$ of the equation (\ref{9}) and (\ref{10}) respectively, we can get the values of
  $p_i$ and $q_{n+i-m}$ are uniquely defined, where $i=0,1,\cdots,{m-1}$. We also can derive the value of $q_j$ is uniquely defined, where $j=1,2,\cdots,{n-m-1}$.

For the value of $q_0$, We  compare the coefficients of $x^{2m+1}$ and  $x^{n+2m+2}$ of the equation (\ref{10}) respectively, we have

  \begin{equation}\label{11}
  \begin{aligned}
2q_{2m+1} b_0+\cdots+ 2q_0 b_{2m+1}=\left({(2m+2)q_{2m+2}(p_0^2-q_0)+\cdots+}\right. \\
\left.{ q_1(2 p_m p_{m+1}-q_{2m+2})}\right) \\
 \end{aligned}
\end{equation}

\begin{equation}\label{12}
 \begin{aligned}
2b_{2m+1} q_{n+1}+\cdots+ 2b_{n}q_{2m+2}=\left({(n+1)q_{n+1}(p_{m+1}^2-q_{2m+2})+\cdots+}\right. \\
\left.{(-2m-2)q_{2m+2}q_{n+1}}\right) \\
\end{aligned}
\end{equation}

the coefficient of $q_{0}$ in the linear equation (\ref{11}) is $2b_{2m+1}+(2m+2)q_{2m+2}$, while $2b_{2m+1}+(2m+2)q_{2m+2}=(n+1)p_{m+1}^2\neq 0$ which  derive from the equation (\ref{12}). Therefore the value of $q_{0}$ is uniquely defined. Finally the polynomial $P(x)$ and $Q(x)$ are determined, we complete the proof of the lemma in the case $n>2m+1$.

If $n<2m+1$, the degree of $P^2(x)-Q(x)$ is smaller than the degree of $P^2(x)$.
Thus, the coefficients of some higher terms of
$P^2(x)$ and $Q(x)$  are same, namely,
\begin{equation}\label{13}
(P^2)^{(n+i)}(x)=(Q)^{(n+i)}(x), \quad 2\leq i\leq 2m+2-n.
\end{equation}
Let us separate $P(x)$ and $Q(x)$ in the form
\begin{displaymath}
P(x)=\sum_{i=0}^{m+1}p_i x^i,\qquad Q(x)=\sum_{i=0}^{2m+2}q_i x^i.
\end{displaymath}
Comparing the coefficients of the highest degree terms of polynomials
on each side of equation (\ref{9}) and (\ref{13}), we have
 \begin{displaymath}
2a_m q_{2m+2}=(4m+4)p_{m+1}q_{2m+2}, \qquad p_{m+1}^2=q_{2m+2}.
\end{displaymath}
Thus  $p_{m+1}=\frac{a_m}{2(m+1)}$ are uniquely defined.

More generally, by comparing the coefficients of $x^{3m+2-i}$ and $x^{2m+2-i}$ of the equation  (\ref{9}) and (\ref{13}) respectively, we have $ p_{m+1-i}=\frac{a_{m-i}}{2(m+1-i)}$ and  the valve of $ q_{2m+2-i}$ is uniquely defined, where $0\leq i\leq {2m-n}$.

Then we compare the coefficients of  $x^{n+m+1}$ and $x^{n+2m+2}$ of the equation (\ref{9}) and  (\ref{10}), we have $ p_{n-m}=\frac{a_{n-m-1}}{2(n-m)}+\frac{(n-2m-1)b_n}{2(n-m)a_m}$, and the valve of $ q_{n+1}$ is uniquely defined. Repeating the above process, we derive that $ p_{n-m-2}£¬\cdots£¬ p_1$ and $q_{n-1}£¬\cdots£¬q_{m+2}$ are uniquely defined.

For the value of $ p_0$, we can compare the coefficient of $x^n$ of equation  (\ref{9}), then we have a linear equation for $ p_0$ while the coefficient of  $ p_0$ can only be $-b_n$ or $\frac{n-4m-3}{m+1}b_n$, therefore the value of $ p_0$ is uniquely defined, then the values of $ q_{m+1},q_m,\cdots,q_0$ which are depend on the value of $ p_0$ are uniquely defined. Finally the polynomial $P(x)$ and $Q(x)$ are determined, we complete the proof of the lemma.

We know from the above discuss, if an $(m,n)$-Lienard system (\ref{1}), where $n\neq2m+1$,
has a hyperelliptic curve
$(y+P(x))^2-Q(x)=0,$
then the system  can only  has this  hyperelliptic curve. Thus, there are at most two points in the  hyperelliptic limit cycles of the system when we fix the value of $x$ which means no  hyperelliptic limit cycle can contained other  hyperelliptic limit cycle. Therefore, the hyperelliptic limit cycles  only have non-nested configuration. (see Fig.2)

\begin{figure}[H]
  \centering
  \includegraphics[width=2.5in]{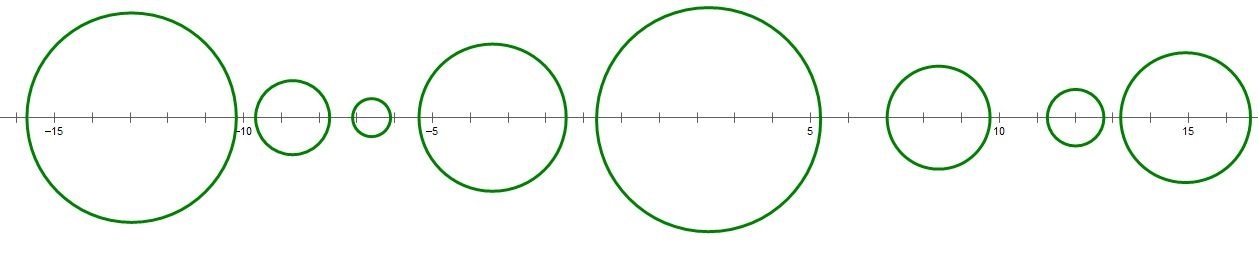}\\
  \label{fig}{Fig.2.The configuration of hyperelliptic limit cycles.}
\end{figure}

\section{The Proof of the Results about Upper Bounds}

By the argument of Lemma \ref{lemma10}, we know a system (\ref{1}) in the case  $n\neq2m+1$
has a hyperelliptic curve $(y+P(x))^2-Q(x)=0$,
then the system  can only  has this  hyperelliptic curve. Take the polynomial $P$, $Q$ of the  hyperelliptic curve in the form
$$P(x)=\prod_{i=1}^a(x-x_i)^{\alpha_i+1}\prod_{j=1}^b(x-y_j)^{\beta_j+1}\prod_{l=1}^c(x-z_l)^{\gamma_l+1}$$
$$Q(x)=\prod_{i=1}^a(x-x_i)\prod_{j=1}^b(x-y_j)^{\omega_j+2},$$
where $a,b,c,\alpha_i,\beta_j,\gamma_l,\omega_j\geq0$, $\alpha=\sum\limits_{i=1}^a{\alpha_i}$, $\beta=\sum\limits_{j=1}^b{\beta_j}$, $\gamma=\sum\limits_{l=1}^c{\gamma_l}$, and $\omega=\sum\limits_{j=1}^b{\omega_j}$. We set $x_i\neq y_j\neq z_l$, $x_1\neq x_2\neq\cdots\neq x_a$, $y_1\neq y_2\neq\cdots\neq y_b$ and $ z_1\neq z_2\neq\cdots\neq z_c$. %for $i=1,2,\cdots,a$, $j=1,2,\cdots,b$, $l=1,2,\cdots,c$, $x_1\neq x_2\neq\cdots\neq x_n$
If $2\beta_i>\omega_i$, $i=1,2,\cdots,b$, then we replace $\beta_i$ and $\omega_i$ with $\beta_i^-$ and $\omega_i^-$ respectively, and use $b^-$ denotes the number of i which satisfied $2\beta_i>\omega_i$.
Otherwise, we replace $\beta_i$ and $\omega_i$ with $\beta_i^+$ and $\omega_i^+$, and use $b^+$ denotes the number of i which satisfied $2\beta_i\leq\omega_i$, then $b=b^++b^-$, $\omega=\omega_i^++\omega_i^-$, $\beta=\beta_i^++\beta_i^-$.

For proving the result of upper bounds, firstly, we discuss the case $m+2\leq n\leq2m-2$.

If
$F(x, y)=(y + P(x))^2-Q(x)=0$ is an invariant algebraic curve of system $(1)$, it is necessary that $P(x)$ has degree $m+1$, and $P^2(x)-Q(x)$ has degree $n+1$, thus
$a+b+c+\alpha+\beta+\gamma=m+1,$
$a+2b+\omega=2m+2$,
and $ln(P^2/Q)=O(x^{n-2m-1})$ , which implies that
 \begin{eqnarray}\label{14}
 \sum\limits_{i=1}^a(2\alpha_i+1)x_i^j+\sum\limits_{i=1}^{b^-}(2\beta_i^--\omega_i^-)y_i^j+\sum\limits_{i=1}^c(2\gamma_i+2)z_i^j
=\sum\limits_{i=1}^{b^+}(\omega_i^+-2\beta_i^+)y_i^j,
 \end{eqnarray}
 where $ j=1,2,\cdots,2m-n$. Assume
 \begin{equation}
f(x)=\prod\limits_{i=1}^a(x-x_i)^{2\alpha_i+1}\prod\limits_{j=1}^{b^-}(x-y_j)^{2\beta_j^--\omega_j^-}\prod\limits_{l=1}^c(x-z_l)^{2\gamma_l+2},
\end{equation}
 \begin{equation}
g(x)=\prod\limits_{j=1}^{b^+}(x-y_j)^{\omega_j^+-2\beta_j^+}.
\end{equation}
We use $k$ denotes the number of the distinct roots of \emph{$f(x)$ }, $t$ denotes the number of the distinct roots of \emph{$g(x)$}, $s$ denotes the degree of $f(x)$, $\tau$ denotes the degree of $f(x)-g(x)$, and $t_0$ denotes the number of the distinct real roots of $Q(x)$. It is easy to see $s=2\alpha+a+2\beta^--\omega^-+2\gamma+2c$, $\tau=n+1-a-2b-2\beta^+-\omega^-$, and $t=b^+$, $k=a+b^-+c$.

 From \emph{\cite{YLJZ}}, we know the discrimination sequence \emph{$(D_1,D_2,\cdots,D_n)$} of \emph{$f(x)$ } satisfied $D_k\neq0, D_{k+1}=D_{k+2}=\cdots=D_s=0$. When $n$ is even, if $t\geq
 \frac{s-\tau+1}{2}$, then $t_0\leq k \leq m+1-b^+\leq\frac{n}{2}$. Otherwise, $t\leq
 \frac{s-\tau-1}{2}$, for the n-degree polynomial $f(x)$, we have

$$S_l=\left|
 \begin{array}{cccc}
 n & \sum\limits_{i=1}^{b^+}(\omega_i^+-2\beta_i^+)y_i & \cdots & \sum\limits_{i=1}^{b^+}(\omega_i^+-2\beta_i^+)y_i^{l-1} \\
\sum\limits_{i=1}^{b^+}(\omega_i^+-2\beta_i^+)y_i &  \sum\limits_{i=1}^{b^+}(\omega_i^+-2\beta_i^+)y_i^2 & \cdots & \sum\limits_{i=1}^{b^+}(\omega_i^+-2\beta_i^+)y_i^{l} \\
 \cdots & \cdots & \cdots & \cdots \\
\sum\limits_{i=1}^{b^+}(\omega_i^+-2\beta_i^+)y_i^{l-1} & \sum\limits_{i=1}^{b^+}(\omega_i^+-2\beta_i^+)y_i^{l} & \cdots & \sum\limits_{i=1}^{b^+}(\omega_i^+-2\beta_i^+)y_i^{2l-2}\\
 \end{array}
\right|,$$ where $0\leq l\leq \frac{s-\tau+1}{2}$.

From Lemma \ref{lemma5}, we have $D_{t+1}=D_{t+2}=\cdots=D_{\frac{s-\tau+1}{2}}=0$ . If $k\leq\frac{s-\tau+1}{2}$, then $t_0\leq k \leq m+1-\frac{n}{2}<\frac{n}{2}$. Otherwise, from Lemma \ref{lemma6}, we have

 $t_0\leq k-([{\frac{\frac{s-\tau+1}{2}-t}{4}}]\times2+[\frac{\frac{s-\tau+1}{2}-t-[{\frac{\frac{s-\tau+1}{2}-t}{4}}]\times4+1}{2}])\times2\leq k-\frac{s}{2}+\frac{\tau}{2}-\frac{1}{2}+t\leq a+b+c+\frac{n+1}{2}-(a+b+c)-\frac{1}{2}\leq\frac{n}{2}$.

 When $n$ is odd, if $t\geq
 \frac{s-\tau}{2}$, then $t_0\leq k\leq m+1-b^+\leq\frac{n+1}{2}$. Otherwise, $t\leq
 \frac{s-\tau-2}{2}$, we have $D_{t+1}=D_{t+2}=\cdots=D_{\frac{s-\tau}{2}}=0$. If $k\leq\frac{s-\tau}{2}$, then $t_0\leq k <\frac{n+1}{2}$. Otherwise, from Lemma \ref{lemma6}, we have

 $t_0\leq k-([{\frac{\frac{s-\tau}{2}-t}{4}}]\times2+[\frac{\frac{s-\tau}{2}-t-[{\frac{\frac{s-\tau}{2}-t}{4}}]\times4+1}{2}])\times2\leq k-\frac{s}{2}+\frac{\tau}{2}+t\leq a+b+c+\frac{n+1}{2}-(a+b+c)\leq\frac{n+1}{2}$.

Since the  system (\ref{1}) can have at most one hyperelliptic limit curves, and the hyperelliptic limit cycle  should intersect the $x$-axis at two different points  $x_1, x_2,$ where  $x_1, x_2,$ are simple root of $Q(x)$, we have $H(m,n)\leq \frac{t_0}{2}$. This completes the proof of case $m+2\leq n\leq2m-2$.

 When $n=2m-1$, if $m$ is odd, we know from the preliminaries, any  root of $Q(x)$ must be   a root of $P(x)$, while the degree of $P(x)$ is $m+1$, $Q(x)$ can have at most $m$ simple roots, then $H(m,n)\leq\frac{m-1}{2}=\left[\frac{n-1}{4}\right]$. For $m$ is even, if $Q(x)$ have $m$ simple roots,  then there are at most  $\frac{m-2}{2}$ intervals which satisfied $Q(x)>0$ in the interval. Otherwise $Q(x)$ can  have at most $m-1$ simple roots, then $H(m,n)\leq\frac{m-2}{2}=\left[\frac{n-1}{4}\right]$. When $n=2m$, the proof  is similar to the case $n=2m-1$, so we omit it.

Recall that we want to prove
$H(m,n)\leq \left[\frac{m}{2}\right]$  when $ n>{2m+1}$.
Since the  system (\ref{1}) can have at most one hyperelliptic limit curves, and $Q(x)$ can  have not more than $m$ simple roots when  $n>2m+1$, we obtain the upper bound of  $H(m,n)$ is
$\left[\frac{m}{2}\right]$ .

\end{document}